\documentclass[oneside, reqno]{amsart}

\usepackage{xypic}
\input xy
\xyoption{all}
\usepackage{epsfig}
\usepackage{amsthm}
\usepackage{amssymb}
\usepackage{amsmath}
\usepackage{amscd}
\usepackage{mathrsfs}
\usepackage{color}
\usepackage{graphicx} 
\usepackage{subcaption,hyperref} 

\hypersetup{
colorlinks=true,
linkcolor=blue,
filecolor=blue,      
urlcolor=blue,
citecolor=blue,
}

\usepackage[parfill]{parskip}
\usepackage{tikz,pgfplots,caption}
\usetikzlibrary{patterns}

\usepackage[T1]{fontenc}
\usepackage{stackengine}
\stackMath

\newcommand{\bg}{\begin{equation}}
\newcommand{\ed}{\end{equation}}
\newcommand{\bga}{\begin{eqnarray}}
\newcommand{\eda}{\end{eqnarray}}
\newcommand{\pf}{\textbf{Proof:\ }}

\def\cbdu{\par{\raggedleft$\Box$\par}}

\newtheorem {Theorem}  {Theorem}

\numberwithin{Theorem}{section}

\newtheorem {Lemma}[Theorem]  {Lemma}

\theoremstyle{definition}
\newtheorem{Definition}[Theorem]{Definition}
\theoremstyle{remark}
\newtheorem{Remark}[Theorem]{\bf Remark}

\expandafter\chardef\csname pre amssym.def
at\endcsname=\the\catcode`\@ \catcode`\@=11
\def\undefine#1{\let#1\undefined}
\def\newsymbol#1#2#3#4#5{\let\next@\relax
 \ifnum#2=\@ne\let\next@\msafam@\else
 \ifnum#2=\tw@\let\next@\msbfam@\fi\fi
 \mathchardef#1="#3\next@#4#5}
\def\mathhexbox@#1#2#3{\relax
 \ifmmode\mathpalette{}{\m@th\mathchar"#1#2#3}%
 \else\leavevmode\hbox{$\m@th\mathchar"#1#2#3$}\fi}
\def\hexnumber@#1{\ifcase#1 0\or 1\or 2\or 3\or 4\or 5\or 6\or 7\or 8\or
 9\or A\or B\or C\or D\or E\or F\fi}

\font\teneufm=eufm10 \font\seveneufm=eufm7 \font\fiveeufm=eufm5
\newfam\eufmfam
\textfont\eufmfam=\teneufm \scriptfont\eufmfam=\seveneufm
\scriptscriptfont\eufmfam=\fiveeufm

\catcode`\@=\csname pre amssym.def at\endcsname

\newcounter{remark}
\newcommand{\R}{\mathbf{R}}

\def  \R   {{\mathbb R}}

\def  \T   {{\mathbb T}}

\def  \12  {{\frac{1}{2}}}

\def\build#1_#2^#3{\mathrel{\mathop{\kern 0pt#1}\limits_{#2}^{#3}}}

\begin{document}

\title[On weak solutions to electron-MHD]{On uniqueness and helicity conservation of weak solutions to the electron-MHD system}


\author [Mimi Dai]{Mimi Dai}

\address{Department of Mathematics, Stat. and Comp.Sci., University of Illinois Chicago, Chicago, IL 60607, USA}
\email{mdai@uic.edu} 

\author [Jacob Krol]{Jacob Krol}

\address{Honors College \& Department of Mathematics, Stat. and Comp.Sci., University of Illinois Chicago, Chicago, IL 60607, USA}
\email{jkrol5@uic.edu} 

\author [Han Liu]{Han Liu}

\address{Department of Mathematics, Stat. and Comp.Sci., University of Illinois Chicago, Chicago, IL 60607, USA}
\email{hliu94@uic.edu}

\thanks{The authors were partially supported by NSF grant
DMS--1815069.}





\begin{abstract} We study the weak solutions to the electron-MHD system and obtain a conditional uniqueness result. In addition, we prove conservation of helicity for weak solutions to the electron-MHD system under a geometric condition.

\bigskip


\hspace{0.02cm}CLASSIFICATION CODE: 35Q35, 35D30, 76D03, 76W05.
\end{abstract}

\maketitle

\section{Introduction}

Considered in this paper is the electron magnetohydrodynamics (EMHD) model arising from plasma physics, written as
\begin{equation}\label{emhd}
\begin{cases}
B_t+d_i\nabla\times((\nabla\times B)\times B)=\mu \Delta B,\\
\nabla\cdot B=0, \ t \in \R^3, \ x \in \R^3 (\text{ or }\T^3).
\end{cases}
\end{equation}
In the above system, the vector valued function $B$ represents the magnetic field whereas the coefficient $d_i$ and $\mu$ stand for the ion inertial length and magnetic resistivity, respectively. System (\ref{emhd}) is a special case of the Hall-MHD system
\begin{equation}\label{HMHD}
\begin{cases}
u_t+u\cdot\nabla u-B\cdot\nabla B+\nabla p= \nu \Delta u,\\
B_t+u\cdot\nabla B-B\cdot\nabla u+d_i\nabla\times((\nabla\times B)\times B)=\mu \Delta B,\\
\nabla \cdot u= 0, \ \nabla\cdot B=0,  \  t \in \R^3, \ x \in \R^3 (\text{ or }\T^3).
\end{cases}
\end{equation} of importance in the studies of a wide range of phenomena and topics in physics, e.g., solar flares, geo-dynamo, aurorae and tokamak.
In system \eqref{HMHD}, the ion flow of the plasma is approximated by an incompressible fluid flow, with $u$ denoting the fluid velocity, $p$ the fluid pressure and $\nu$ the viscosity coefficient. System \eqref{HMHD} differs from standard MHD systems by the term $d_i\nabla\times((\nabla\times B)\times B)$ describing the Hall effect, which becomes significant at sub-ion scale.  It is believed that in this setting the Hall effect alters Alfv\'en's ``frozen-in'' theorem for the standard MHD, a violation of which is essential to the magnetic reconnection process, i.e., the topological reorganization of magnetic field lines, widely observed in space plasmas. At scales $\ell \ll d_i,$ system \eqref{HMHD} reduces to system \eqref{emhd}, as the ions and electrons become decoupled, causing the magnetic field lines to be frozen into the electronic fluid only. For more physical background, we refer readers to \cite{G1, L, MG}.

It is sometimes assumed that there is a magnetic vector potential $A,$ satisfying $B=\nabla \times A,$ which, under the assumption of Coulomb gauge, can be chosen such that $\nabla\cdot A=0.$ Thus, $A$ can be recovered from $B$ through Biot-Savart law, i.e., $$A=\nabla\times (-\Delta)^{-1}B.$$ Formally, $A$ satisfies the following system of equations
\begin{equation}\label{EAMHD}
\begin{cases}
A_t - d_i (\nabla \times B)\times B = \mu \Delta A,\\
\nabla \times A = B,\\
\nabla \cdot B =0, \ t \in \R^+, \ x \in \R^3 (\text{ or }\T^3).
\end{cases}
\end{equation}
In our paper, we shall work with the above form of the EMHD system as well.

There is a sizable literature on the mathematical studies of Hall-MHD and EMHD systems. Global existence of weak solutions was established in \cite{ADFL, CDL, DS}, while several well-posedness results can be found in \cite{CDL, CWW, D2, D3, DT, DL2, KL}. In \cite{CW, JO}, ill-posedness results were obtained whereas non-uniqueness of weak solutions was proven in \cite{D4}. In addition, the asymptotic behavior of solutions was studied in \cite{CS, DL}. For various regularity and blow-up criteria, readers are referred to \cite{CL, CW, D1, FFNZ, FLN, HAHZ, WZ1, WL, Y3, Y4, YZ0, Z}.

A fundamental result is the global existence of Leray-Hopf type weak solution to system \eqref{emhd}, which can be proven via a standard Galerkin approximation procedure (cf. \cite{CDL}). The Leray-Hopf type weak solution is defined as follows.
\begin{Definition}\label{weak-sol}
$B$ is said to be a weak solution of (\ref{emhd}) on $[0, T]$ if $B$ is divergence-free in the sense of distributions and satisfies following integral equation
\begin{equation}\notag
\int_0^T\int_{\mathbb R^3} \Big(B\cdot \varphi_t+(B\otimes B):\nabla\nabla\times \varphi \Big) \,\mathrm dx\,\mathrm dt=\mu \int_0^T\int_{\mathbb R^3}\nabla B:\nabla \varphi\,\mathrm dx\,\mathrm dt
\end{equation}
for any $\varphi\in \mathscr{D}([0,T]\times \mathbb R^3).$

Moreover, a weak solution $B$ is called a Leray-Hopf type solution if $$B \in L^\infty(0,T; L^2(\mathbb R^3))\cap L^2(0,T; H^1(\mathbb R^3))$$
and the energy inequality
\begin{equation}\label{energy1}
\|B(t)\|_{L^2}^2+2\mu \int_{t_0}^t\|\nabla B\|_{L^2}^2 \mathrm dt\leq \|B_0\|_{L^2}^2
\end{equation}
holds for almost every $t_0\in[0,T]$ and $t \in (t_0, T].$
\end{Definition}
The uniqueness of Leray-Hopf type solutions, however, remains an open question. In fact, on the negative side, non-uniqueness of weak solutions to system \eqref{emhd} in the Leray-Hopf class $L^\infty(0,T; L^2(\mathbb R^3))\cap L^2(0,T; H^1(\mathbb R^3))$ has been proven in \cite{D4} via the celebrated convex integration method. 

The first result of this paper concerns the positive side of the uniqueness question for the Leray-Hopf type solutions. It is a so-called weak-strong uniqueness result, stated as follows.
\begin{Theorem}\label{thm1}
Given two divergence-free vector fields $B_0 ^{1}, B_0^{2}\in L^2(\mathbb R^3),$ denote by $B^1(t)$ and $B^2(t)$ two Leray-Hopf type weak solutions to (\ref{emhd}) on $[0,T)$ generated by $B^1_0$ and $B^2_0,$ respectively. If
$$\nabla\times B^1\in L^q(0,T; B^r_{p,\infty})$$
with 
$$\frac2q+\frac3p=1+r, \ \ \frac{3}{1+r}<p\leq \infty, \ r\in(0,1], \  (p,q)\neq (\infty, 1), $$
then for $t \in (0,T),$ the inequality
\begin{equation}\notag
\|B^1(t)-B^2(t)\|_{L^2}^2\leq \|B_0^1-B_0^2\|_{L^2}^2\exp\Big\{C\big(t+\|\nabla\times B^1\|_{L^q(0,t; B^r_{p,\infty})}\big)\Big\}
\end{equation}
holds. In particular, $B^1=B^2$ a.e. on $[0,T)\times \mathbb R^3$ provided that $B_0^1=B_0^2$. 
\end{Theorem}
\begin{Remark}
We note that system \eqref{emhd} is invariant under the following scaling transformation $$B(x, t) \mapsto B_\lambda(x,t):= B(\lambda x, \lambda^2 t).$$ The  condition $\nabla\times B \in L^q(0,T; B^r_{p,\infty})$ is consistent with the above scaling symmetry. 
\end{Remark}
\begin{Remark}
In \cite{CL}, it was shown that a weak solution to system \eqref{emhd} is regular, thus unique, on $[0,T]$ if and only if 
$$\nabla B \in L^q(0,T; L^p(\mathbb R^3)) \text{ for }\frac{2}{q}+\frac{3}{p} \leq 1.$$ This regularity criterion is consistent with our conditional uniqueness result.  
\end{Remark}

As noted before, the Hall-MHD system is an essential model in interpreting the magnetic reconnection process, responsible for celestial events from aurorae caused by magnetic substorms in planetary magnetospheres to the violent solar flares. Since magnetic reconnection features topological changes of magnetic field lines, it is therefore of interest to study the magnetic helicity 
$$\mathcal{H}(t) = \int_{\R^3} \big(A\cdot B\big)(t,x)\mathrm{d}x,$$
which is regarded as a tool to quantify the magnetic topology, i.e., self-linkage and knottedness of magnetic field lines. Clearly, $\mathcal{H}(t)$ is dissipated by the diffusive term $\mu \Delta B$ in the resistive setting. Yet, besides the presence of magnetic resistance, the lack of regularity of the solution can also cause the dissipation of $\mathcal{H}(t),$ known as anomalous dissipation. The concept was first postulated by Onsager in the context of hydrodynamics and has been validated for the Navier-Stokes and Euler equations (cf. \cite{I}).

Our second result addresses the issue of magnetic helicity conservation, i.e., preservation of the magnetic topology, for weak solutions to system \eqref{EAMHD}. More specifically, we shall give a set of conditions on the weak solutions to system \eqref{EAMHD} so that for $\phi \in \mathscr{D}([0,T] \times \R^3)$ and $t \in (0,T],$ the following generalized helicity equality, which implies the absence of anomalous dissipation, holds
\begin{equation}\notag
\begin{split}
\int_{\R^3 \times \{t\}} A \cdot B \phi + 2\mu \int^t_0 \int_{\R^3} \nabla A : \nabla B \phi = & \int_{\R^3 \times \{0\}} A \cdot B \phi+\int^t_0 \int_{\R^3} A \cdot B (\phi_t + \mu \Delta \phi)\\
& + d_i\int^t_0 \int_{\R^3} ((\nabla \times B) \times B) \cdot (\nabla \phi \times A).
\end{split}
\end{equation}
Our result is as follows.
\begin{Theorem}\label{thm2}
Let $s \in C^{\frac{1}{2}}([0,T] \times \R^3)$ and $(A, B)$ be a weak solution to \eqref{EAMHD} satisfying 
\begin{align}
&A \in C_w H^1 \cap L^3 W^{1,\frac{9}{2}} \cap L^2 H^2;\label{acond}\\
&\nabla(\nabla \times A) \in \big(L^3 L^{9/5}\cap L^{3/2}L^{18/5}\big)((0,T)\times\R^3/\text{Graph}(s))_{\text{loc}},\label{bcond}
\end{align}
then $(A, B)$ satisfies the general helicity identity.
\end{Theorem}

\begin{Remark}
In \cite{S}, it was shown that if a Leray-Hopf weak solution $u$ to the Navier-Stokes equations satisfies $u \in L^3 L^{9/2}$ and $\nabla u$ belongs locally to $L^3 L^{9/5}$ outside a $C^{1/2}$- curve, then for $u$ the generalized energy equality holds. In this paper, we adapt the idea therein to system \eqref{EAMHD}.
\end{Remark}

\begin{Remark}
The spaces $L^3 W^{1, \frac{9}{2}}$ for $A$ and $L^3 L^{9/5}$ for $\nabla B$ are Onsager critical. Due to the asymmetry of the Hall term, we need the additional assumption that $\nabla B \in L^{3/2} L^{18/5}((0,T)\times\R^3/\text{Graph}(s))_{\text{loc}}.$  
\end{Remark}

In the case of the non-resistive EMHD system, conservation of magnetic helicity for weak solutions in the Onsager critical Chemin-Lerner space $\widetilde L^3(0,T; B^{1/3}_{3, c(\mathbb{N})})$ was proven in \cite{DS}. In the appendix of this paper, we shall give a proof of the following variant of the result via Littlewood-Paley theory.
\begin{Theorem}\label{consh}
Let $B \in L^3(0,T; B^{1/3}_{3, c(\mathbb{N})})\cap C_w(0,T; H^{-\frac{1}{2}})$ be a weak solution to the non-resistive EMHD system, then $B$ conserves the magnetic helicity $\mathcal{H}.$
\end{Theorem}

\bigskip

\section{Preliminaries}

\subsection{Notation}
For simplicity, we denote by $L^p X$ the space $L^p(0,T; X(\R^3)),$ where $X$ is a Banach space, and by $L^p L^q((0,T)\times\R^3/\text{Graph}(s))_{\text{loc}}$ the space $$\big\{f \in \mathscr{D}': f\phi \in L^p L^q, \forall \phi \in \mathscr{D}((0,T)\times\R^3/\text{Graph}(s))\big\},$$ where $s \in C^{1/2}([0,T] \times \R^3).$ 

For shortness, we sometimes write $\|\cdot \|_{L^p}$ as $\|\cdot \|_p$. For two matrices $X, Y \in \mathcal{M}_{3 \times 3},$ the notation $X:Y$ refers to $\text{Tr}[X \otimes Y].$

\subsection{Vector calculus identities}\label{vecid}

Let $A$ and $B$ be vector valued functions, and $\varphi$ be a scalar function. We shall use the following identities --
\begin{gather*}
\nabla(\varphi A)= \nabla \varphi \otimes A+ \varphi \nabla A;\\
\nabla \cdot (\varphi A)=\nabla \varphi \cdot A+ \varphi \nabla \cdot A;\\
\nabla\times (\varphi A)= \varphi(\nabla\times A)+(\nabla \varphi)\times A;\\
\nabla \times (A \times B) = A (\nabla \cdot B) - B(\nabla \cdot A)+ (B \cdot \nabla) A - (A \cdot \nabla) B;\\
(\nabla\times A)\times B = A \times (\nabla \times B) + (A \cdot \nabla) B + (B \cdot \nabla)A - \nabla (A \cdot B) .
\end{gather*}
In particular, setting $A=B$ in the last inequality above yields
$$(\nabla \times B) \times B = \nabla \cdot (B \otimes B) -\frac{1}{2}\nabla |B|^2.$$
We also use the facts that $\nabla \times (\nabla B) =0$ and $(A \times B) \cdot A =0.$

\subsection{Besov spaces via Littlewood-Paley theory}

For $s \in \R$ and $1 \leq p,q \leq \infty,$ we define the inhomogeneous Besov space $ B^s_{p,q}$ as 
$$\dot B^s_{p,q}(\R^n) = \big\{f \in \mathscr{S}'(\R^n) : \|f\|_{B^{s}_{p,q}(\R^n)} < \infty \big \},$$ with the norm given by
\begin{equation}\notag
\|f\|_{B^{s}_{p,q}(\R^n)}=
\begin{cases}
\displaystyle\big(\sum_{j \geq -1 } (2^{sj}\|\Delta_j f \|_{L^p(\R^n)})^q \big)^{\frac{1}{q}}, \ \text{ if }1 \leq q < \infty,\\
\displaystyle \sup_{j \geq -1} (2^{sj}\|\Delta_j f \|_{L^p(\R^n)}), \ \text{ if } q= \infty,
\end{cases}
\end{equation}
where $\Delta_j$ is the $j$-th inhomogeneous Littlewood-Paley projection.

In Theorem \ref{consh}, the notation $B^s_{p, c(\mathbb{N})}$ refers to the union of all Besov spaces $B^s_{p, q}$ with finite $q,$ endowed with the norm $\|\cdot\|_{B^s_{p, \infty}}.$ 

We shall attach a brief review of Littlewood-Paley theory in the appendix.

\subsection{A decomposition lemma}

The following lemma, found in \cite{CMZ}, turns out to be useful in the proof of Theorem \ref{thm1}.
\begin{Lemma}\label{le-decom}
Assume $B \in L^q(0,T; B^r_{p,\infty})$ with $\frac2q+\frac3p=1+r$, $\frac3{1+r}<p\leq \infty$, $r\in(0,1]$, and $(p,q)\neq (\infty, 1)$. Then $B$ can be decomposed as $B=B^\ell+B^h$ with
$$B^\ell\in L^1(0,T; \mathrm {Lip})\ \mathrm {and}\ B^h\in L^{q'}(0,T; L^{p'})$$
for some $p'$ and $q'$ satisfying $\frac2{q'}+\frac3{p'}=1$, $p'>3$. Moreover, for $t \in (0, T],$ the following estimate holds --
\begin{equation}\notag
\int^t_0 \left(\|\nabla B^\ell(\tau)\|_\infty+\|B^h (\tau)\|_{p'}^{q'}\right)  \mathrm{d}\tau \leq   C\int_0^t \left(1+\|\nabla\times B^1(\tau)\|_{B^r_{p,\infty}} \right)^q \mathrm{d}\tau.
\end{equation}
\end{Lemma}

\bigskip

\section{Proof of Theorem \ref{thm1}}

\pf We assume that $B^1(t)$ and $B^2(t)$ are two Leray-Hopf type solutions to system (\ref{emhd}) with initial data $B^1_0$ and $B^2_0,$ respectively and denote $$Z(t):=B^2(t)-B^1(t).$$ Taking the inner product of $B^1$ equation with $B^2$ and vice versa, then  integrating over $\mathbb R^3\times [0,t]$ yields the following equality. (This procedure can be done rigorously using Galerkin approximations.)
\begin{equation}\notag
\begin{split}
&\int_{\mathbb \R^3\times \{t\}}B^1\cdot B^2 -\int_{\mathbb \R^3 \times \{0\}}B^1 \cdot B^2\\
=&\ d_i\int_0^t\int_{\mathbb R^3}\Big(\nabla\times \big((\nabla\times B^1\big)\times B^1)\cdot B^2+\nabla\times \big((\nabla\times B^2\big)\times B^2)\cdot B^1\Big)\\
&+ \mu \int_0^t\int_{\mathbb R^3} (\Delta B^1\cdot B^2+\Delta B^2\cdot B^1)\\
=&: I_1+I_2.
\end{split}
\end{equation}
Integration by parts leads to
$$ I_2=-2\mu\int_0^t\int_{\mathbb R^3}\nabla B^1\cdot \nabla B^2.$$

Integrating by parts and using vector identities from Section \ref{vecid}, we can rewrite $I_1$ as follows.
\begin{equation}\notag
\begin{split}
I_1=& d_i \int_0^t\int_{\mathbb R^3}\Big(\nabla\times \big((\nabla\times B^1)\times B^1\big)\cdot Z+\nabla\times \big((\nabla\times B^1)\times B^1\big)\cdot B^1\Big) \\
&+ d_i \int_0^t\int_{\mathbb R^3}\nabla\times \Big(\big((\nabla\times B^2)\times B^2 \big)\cdot (-Z)+\nabla\times \big((\nabla\times B^2)\times B^2\big)\cdot B^2\Big)\\
=& d_i \int_0^t\int_{\mathbb R^3}\nabla\times \Big(\big((\nabla\times B^1)\times B^1\big)\cdot Z-\nabla\times \big((\nabla\times B^2)\times B^2\big)\cdot Z \Big)\\
=& - d_i \int_0^t\int_{\mathbb R^3}\Big(\nabla\times \big((\nabla\times B^1)\times Z\big)\cdot Z+\nabla\times \big((\nabla\times Z)\times B^2\big)\cdot Z\Big)\\
=& - d_i \int_0^t\int_{\mathbb R^3}\nabla\times ((\nabla\times B^1)\times Z)\cdot Z.
\end{split}
\end{equation}

Summarizing the analysis above provides 
\begin{equation}\label{en1}
\begin{split}
&\int_{\mathbb R^3 \times \{t\}}B^1 \cdot B^2-\int_{\mathbb R^3 \times \{ 0\}}B^1 \cdot B^2\\
=& -2 \mu \int_0^t\int_{\mathbb R^3}\nabla B^1: \nabla B^2+d_i \int_0^t\int_{\mathbb R^3}\nabla\times \big((\nabla\times B^1)\times Z \big)\cdot Z.
\end{split}
\end{equation}

Since $B^1(t)$ and $B^2(t)$ are Leray-Hopf type solutions, they satisfy the following energy inequalities
\begin{equation}\label{en2}
\begin{split}
\|B^1(t)\|_{2}^2+2\mu \int_0^t\|\nabla B^1(\tau)\|_{2}^2 \mathrm{d}\tau \leq \|B_0^1\|_{2}^2, \\
\|B^2(t)\|_{2}^2+2\mu \int_0^t\|\nabla B^2(\tau)\|_{2}^2 \mathrm{d}\tau \leq \|B_0^2\|_{2}^2.
\end{split}
\end{equation}
In view of (\ref{en1}) and (\ref{en2}), we can derive the following energy inequality for $Z(t).$
\begin{equation}\label{en3}
\begin{split}
& \|Z(t)\|_{2}^2+2\mu \int_0^t\|\nabla Z(\tau)\|_{2}^2\mathrm{d}\tau\\
=&\ \|B^1(t)\|_{2}^2+\|B^2(t)\|_{2}^2 +2 \mu \int_0^t\Big(\|\nabla B^1(\tau)\|_{2}^2+ \|\nabla B^2(\tau)\|_{2}^2\Big)\mathrm{d}\tau \\ &-2\int_{\mathbb R^3 \times \{t\}}B^1 \cdot B^2 -4\mu \int_0^t \int_{\mathbb R^3}\nabla B^1 : \nabla B^2\\\
\leq &\ \|B^1_0 -B^2_0\|_{2}^2-2 d_i \int_0^t\int_{\mathbb R^3}\nabla\times ((\nabla\times B^1)\times Z)\cdot Z.
\end{split}
\end{equation}

Owing to the vector calculus identity
\begin{equation}\notag
\begin{split}
\nabla\times ((\nabla\times B^1)\times Z)=&\ \nabla\times B^1(\nabla\cdot Z)-Z(\nabla\cdot\nabla\times B^1)\\
&+Z\cdot\nabla(\nabla\times B^1)-\nabla\times B^1\cdot\nabla Z\\
=&\ Z\cdot\nabla(\nabla\times B^1)-\nabla\times B^1\cdot\nabla Z
\end{split}
\end{equation}
we can write the flux term in (\ref{en3}) as
\begin{equation}\notag
\begin{split}
\int_0^t\int_{\mathbb R^3}\nabla\times ((\nabla\times B^1)\times Z)\cdot Z =&\ \int_0^t\int_{\mathbb R^3}\Big((Z\cdot\nabla)(\nabla\times B^1) -(\nabla\times B^1\cdot\nabla)Z\Big)\cdot Z \\
=&\ \int_0^t\int_{\mathbb R^3}(Z\cdot\nabla)(\nabla\times B^1)\cdot Z.
\end{split}
\end{equation}

By Lemma \ref{le-decom}, $\nabla\times B^1\in L^q(0,T; B^r_{p, \infty})$ can be decomposed as
$$\nabla\times B^1=\mathcal B^\ell+\mathcal B^h,$$
where $\mathcal B^\ell\in L^1(0,T; \mathrm {Lip})$ and $\mathcal B^h\in L^{q'}(0,T; L^{p'})$ for some $p' >3$ and $q'$ satisfying $\frac2{q'}+\frac3{p'}=1.$ Therefore, the flux term can be written as 
\begin{equation}\notag
\int_0^t\int_{\mathbb R^3}(Z\cdot\nabla)(\nabla\times B^1)\cdot Z =\ \int_0^t\int_{\mathbb R^3}(Z\cdot\nabla)\mathcal B^\ell\cdot Z + \int_0^t\int_{\mathbb R^3}(Z\cdot\nabla)\mathcal B^h\cdot Z.
\end{equation}

The estimate for the first integral on the right hand side is given by
\begin{equation}\label{en4}
\begin{split}
\left| \int_0^t\int_{\mathbb T^3}(Z\cdot\nabla)\mathcal B^\ell\cdot Z\right|\leq \int_0^t\|Z(\tau)\|_{2}^2\|\nabla \mathcal B^\ell(\tau)\|_{\infty}\mathrm{d}\tau.
\end{split}
\end{equation}
To estimate the second integral on the right hand side, we  integrate by parts and apply H\"older's inequality, Gagliardo-Nirenberg inequality and Young's inequality.
\begin{equation}\label{en5}
\begin{split}
\left|\int_0^t\int_{\mathbb R^3}(Z\cdot\nabla)\mathcal B^h\cdot Z \right|=&\left| \int_0^t\int_{\mathbb R^3}(Z\cdot\nabla)Z\cdot \mathcal B^h\right|\\
\leq &\  \int_0^t\|\nabla Z(\tau)\|_{2}\|Z(\tau)\|_{\frac{2p'}{p'-2}}\|\mathcal B^h (\tau)\|_{p'} \mathrm{d} \tau \\
\leq &\ C\int_0^t\|\nabla Z(\tau)\|_{2}^{1+\frac3{p'}}\|Z(\tau)\|_{2}^{1-\frac{3}{p'}}\|\mathcal B^h(\tau)\|_{{p'}}\mathrm{d} \tau\\
\leq &\ C\left(\int_0^t\|Z(\tau)\|_{2}^2\|\mathcal B^h(\tau)\|_{p'}^{q'}\mathrm{d}\tau \right)^{\frac1{q'}}\left(\int_0^t\|\nabla Z (\tau)\|_{2}^2 \mathrm{d}\tau\right)^{1-\frac{1}{q'}}\\
\leq &\ C\int_0^t\|Z(\tau)\|_2^2\|\mathcal B^h(\tau)\|_{p'}^{q'}\mathrm{d}\tau +\frac{\mu}{2 d_i}\int_0^t\|\nabla Z(\tau)\|_{2}^2 \mathrm{d}\tau.
\end{split}
\end{equation}

Combining (\ref{en3})-(\ref{en5}) and invoking Lemma \ref{le-decom} yield
\begin{equation}\notag
\begin{split}
& \|Z(t)\|_{L^2}^2+\mu \int_0^t\|\nabla Z(\tau)\|_{L^2}^2\mathrm{d} \tau \\ \leq &\ \|B^1_0-B^2_0\|_{2}^2+C\int_0^t\|Z(\tau)\|_{2}^2\left(\|\nabla \mathcal B^\ell(\tau)\|_{\infty}+\|\mathcal B^h(\tau)\|_{p'}^{q'}\right)\mathrm{d}\tau\\
\leq &\ \|B^1_0-B^2_0\|_{2}^2 +C\int_0^t\|Z(\tau)\|_{L^2}^2\left(1+\|\nabla\times B^1(\tau)\|_{B^r_{p,\infty}}\right)^q\mathrm{d}\tau.
\end{split}
\end{equation}
By Gr\"onwall's inequality, we have  
\begin{equation}\notag
\|Z(t)\|_{L^2}^2\leq \|B^1_0-B^2_0\|_{L^2}^2\exp \left\{C\int_0^t \left(1+\|\nabla\times B^1(\tau)\|_{B^r_{p,\infty}} \right)^q \mathrm{d}\tau \right\}.
\end{equation}
Therefore, it follows that $Z(t)=B^2(t)-B^1(t)=0$ for $t\in[0,T)$ if $B_0^1=B^2_0.$
\cbdu

\bigskip

\section{Proof of Theorem \ref{thm2}}

\pf We start the proof, which is based on an approximation argument, by fixing a mollifier $\eta \in C^\infty_0(B(0,1))$ such that $\eta \geq 0$ and $\int \eta =1.$ For a vector field $B \in (\mathscr{D}'(\R^3))^3,$ we denote
$$B_\delta(x):= \delta^{-3}\int_{\R^3}\eta(\delta^{-1}y) B(x-y)\mathrm{d}y.$$

We define the extension of $s$ as
\begin{equation}\notag
s^{\text{ext}}(t)=
\begin{cases}
s(0), t <0,\\
s(t), 0 \leq t \leq T,\\
s(T), t >T.
\end{cases}
\end{equation}
Clearly, $s^{\text{ext}} \in C^{\frac{1}{2}}(\R \times \R^3).$ Let $\eta$ be a mollifier. We approximate the graph of $s$ by
$$s_\varepsilon = \varepsilon^{-2}\int_{\R} \eta(\varepsilon^{-2} \tau) s^{\text{ext}}(t-\tau)\mathrm{d}\tau,$$
which satisfies the following inequalities --
\begin{equation}\label{ineq'}
\sup_{0 \leq t \leq T}|s(t)-s_\varepsilon (t) | <\varepsilon, \ \sup_{0 \leq t \leq T}|s'_\varepsilon (t) | \leq \varepsilon^{-1}.
\end{equation}
To cut off the graph of $s,$ we let $\chi \in C^\infty(\R^3)$ be such that $\chi \equiv 1$ on $\R^3/B(0,3),$ $0 \leq \chi \leq 1$ on $B(0,3)/B(0,2)$ and $\chi \equiv 0$ on $B(0,2)$ and set
$$\chi_\varepsilon(t,x):= \chi\Big(\frac{x-s_\varepsilon(t)}{\varepsilon}\Big).$$ 
We will use the fact that $$\text{supp }\nabla \chi_\varepsilon \subset \big\{(t,x) : |x-s(t)| \leq 3 \varepsilon \big\}$$ and the following bound on the derivatives of $\chi_\varepsilon$ --
\begin{equation}\label{ineqd}
\sup_{0\leq t \leq T}\|D^\gamma_x \chi_\varepsilon \|_p  \sim \varepsilon^{\frac{3}{p}-\gamma}.
\end{equation}

For a weak solution $(A,B)$ to system \eqref{EAMHD} and $\varphi, \psi \in (\mathscr{D}([0,T] \times \R^3))^3,$ we have 
\begin{equation}\notag
\begin{split}
\int_{\R^3\times \{t\}} A \cdot \varphi  -\int^t_0\int_{\R^3} A \cdot \varphi_t  = & \int_{\R^3 \times \{0\}} A \cdot \varphi - d_i \int^t_0 \int_{\R^3} (B \otimes B): \nabla \varphi\\
&+ \frac{d_i}{2} \int^t_0 \int_{\R^3} |B|^2 \nabla \cdot \varphi - \mu \int^t_0 \int_{\R^3} \nabla A : \nabla \varphi,
\end{split}
\end{equation}
\begin{equation}\notag
\begin{split}
\int_{\R^3\times \{t\}} B \cdot \psi  -\int^t_0\int_{\R^3} B \cdot \psi_t  = & \int_{\R^3 \times \{0\}} B \cdot \psi + d_i\int^t_0  \int_{\R^3} ((\nabla \times B) \times B):(\nabla \times \psi)\\
& - \mu \int^t_0 \int_{\R^3} \nabla B : \nabla \psi.
\end{split}
\end{equation}

Choosing $\varphi= (B_\delta\phi \chi_\varepsilon)_\delta $ and $\varphi= (A_\delta\phi \chi_\varepsilon)_\delta$ with $\phi \in \mathscr{D}([0,T]\times \R^3)$ and summing the above identities, we obtain an identity of the form
$$I- II= III -IV +V +VI -VII,$$
where the terms $I$ -- $VII$ are defined as
\begin{gather*}
I:= 2 \int_{\R^3 \times \{t\}} A_\delta \cdot B_\delta \phi \chi_\varepsilon, \ \ \ II:= \int^t_0 \int_{\R^3} A_\delta \cdot \big(B_\delta  \phi \chi_\varepsilon \big)_t + \big( A_\delta \phi \chi_\varepsilon \big)_t \cdot B_\delta,\\
III := 2 \int_{\R^3 \times \{0\}} A_\delta \cdot B_\delta \phi \chi_\varepsilon, \ \ \ IV:= d_i\int^t_0 \int_{\R^3} (B \otimes B)_\delta : \nabla (B_\delta \phi\chi_\varepsilon),\\
V := d_i \int^t_0 \int_{\R^3} \big( (\nabla \times B) \times B \big)_\delta \cdot \nabla \times (A_\delta \phi\chi_\varepsilon) ,\\
VI:= \frac{d_i}{2}\int^t_0 \int_{\R^3} (B \cdot B)_\delta \nabla \cdot B_\delta \phi \chi_\varepsilon, \\
VII := \mu \int^t_0 \int_{\R^3} \nabla A_\delta : \nabla (B_\delta \phi \chi_\varepsilon)+\nabla(A_\delta\phi \chi_\varepsilon) : \nabla B_\delta.
\end{gather*} 
Our goal is to show that as $\delta, \varepsilon \to 0,$ the above identity converges to the generalized helicity identity. 

To treat terms $I$ and $II$, we exploit the cancellations by integrating $II$ by parts. As a result, the left hand side becomes
\begin{equation}\notag
\begin{split}
\int_{\R^3 \times \{t\}} A_\delta \cdot B_\delta \phi\chi_\varepsilon + \int_{\R^3 \times \{0\}} A_\delta \cdot B_\delta \phi\chi_\varepsilon - \int^t_0 \int_{\R^3} A_\delta \cdot B_\delta (\phi \chi_\varepsilon)_t.
\end{split}
\end{equation} 
Let $\delta, \varepsilon \to 0.$ The first two terms above converge to
$$\int_{\R^3 \times \{t\}} A \cdot B \phi + \int_{\R^3 \times \{0\}} A \cdot B \phi,$$ while for the third term, we write 
\begin{equation}\notag
\begin{split}
\int^t_0 \int_{\R^3} A_\delta \cdot B_\delta (\phi \chi_\varepsilon)_t =& \int^t_0 \int_{\R^3} A_\delta \cdot B_\delta \phi_t \chi_\varepsilon + \int^t_0 \int_{\R^3} A_\delta \cdot B_\delta \phi \ s_\varepsilon'(t) \nabla \chi_\varepsilon\\
=: & II_1 + II_2.
\end{split}
\end{equation}
As $\delta, \varepsilon \to 0,$ $II_1$ converges to its natural limit $$\iint A \cdot B \phi_t.$$

On the other hand, using H\"older's inequality, inequalities \eqref{ineq'} and \eqref{ineqd} along with the fact that $A \in L^\infty L^6$ and $B \in L^2 L^6,$ we estimate $II_2$ as follows.
\begin{equation}\notag
\begin{split}
|II_2| \leq & \|A_\delta\|_{L^\infty L^6}\int^t_0\Bigg(\int_{|x-s_\varepsilon(t)| \leq 3 \varepsilon} |B_\delta|^6|\phi|^3 \ \mathrm{d}x\Bigg)^\frac{1}{3} |s_\varepsilon'(t)|\varepsilon \ \mathrm{d}t\\
\leq & \|A_\delta\|_{L^\infty L^6} \int^t_0\Bigg(\int_{|x-s_\varepsilon(t)| \leq 3 \varepsilon} |B_\delta|^6 \ \mathrm{d}x\Bigg)^\frac{1}{3}\mathrm{d}t\\
\leq & \|A_\delta \|_{L^\infty L^6} \int^t_0\Bigg(\int_{|x-s_\varepsilon(t)| \leq 3 \varepsilon+\delta} |B|^6 \ \mathrm{d}x\Bigg)^\frac{1}{3}\mathrm{d}t,
\end{split}
\end{equation}
which vanishes when $\delta, \varepsilon \to 0.$ 

It is clear that $$III \xrightarrow{\delta, \varepsilon \to 0} 2\int_{\R^3 \times \{0\}} A \cdot B \phi.$$

Introducing the following bilinear form $$r_\delta(B, B) = \delta^{-3} \int_{\R^3} \eta(\delta^{-1}y) \big(B(x -y) -B(x)\big) \otimes \big(B(x -y) -B(x)\big) \mathrm{d}y,$$
we split $IV$ into three parts --
\begin{equation}\notag
\begin{split}
IV = & d_i\iint (B \otimes B)_\delta \nabla (B_\delta \phi \chi_\varepsilon)\\
 = & d_i\iint r_\delta (B,B) \nabla (B_\delta \phi \chi_\varepsilon)+ d_i\iint (B-B_\delta)\otimes(B-B_\delta) \nabla (B_\delta \phi \chi_\varepsilon)\\
& + d_i\iint B_\delta \otimes B_\delta \nabla (B_\delta \phi \chi_\varepsilon)\\
=:& IV_1 + IV_2 + IV_3.
\end{split}
\end{equation}

We note that
$$\|r_\delta (B ,B)\|_{9/4} \leq \delta^{-3}\int \eta(y/\delta)\|B(\cdot )- B(\cdot ) \|_{9/2}^2\ \mathrm{d}y =: R(t, \delta),$$
and the fact that $B \in L^3 L^{9/2}$ implies that $R(t, \delta)$ satisfies the estimate
$$\int^t_0 R^{3/2}(t, \delta) \mathrm{d}t \leq \int^t_0 \int_{\R^3} \delta^{-3} \eta(\delta^{-1}y) \|B(\cdot -y) -B(\cdot)\|_{9/2}^3 \xrightarrow{\delta \to 0} 0,$$
which, along with the condition \eqref{bcond}, yields
\begin{equation}\notag
\begin{split}
|IV_1| \leq & d_i\int^t_0 \|r_\delta(B, B)\|_{9/4} \|\nabla (B_\delta \phi \chi_\varepsilon)\|_{9/5} \\
\leq & \bigg(\int^t_0 R^{3/2}(t, \delta) \mathrm{d}t \bigg)^{2/3} \|\nabla (B_\delta \phi \chi_\varepsilon)\|_{L^{3} L^{9/5}} \xrightarrow{\delta, \varepsilon \to 0} 0.
\end{split}
\end{equation}
By the same argument, we can see that $IV_2$ vanishes, since
\begin{equation}\notag
\begin{split}
|IV_2| \leq & d_i\int^t_0 \|B-B_\delta\|_{9/2}^2\|\nabla (B_\delta \phi \chi_\varepsilon)\|_{9/5} \xrightarrow{\delta \to 0} 0.
\end{split}
\end{equation}

We write $IV_3$ as 
\begin{equation}\notag
\begin{split}
IV_3 = & \frac{d_i}{2}\iint (B_\delta \otimes B_\delta) (B_\delta \nabla \phi \chi_\varepsilon) + \frac{d_i}{2}\iint (B_\delta \otimes B_\delta) (B_\delta  \phi \nabla \chi_\varepsilon)\\
=& IV_{31}+IV_{32},
\end{split}
\end{equation}
where $$ IV_{31} \xrightarrow{\delta, \varepsilon \to 0} \frac{d_i}{2}\iint (B \otimes B):(B \otimes \nabla \phi).$$ 

To estimate $IV_{32},$ we use the estimate \eqref{ineqd} and condition \eqref{acond}. It follows that
\begin{equation}\notag
\begin{split}
|IV_{32}| \leq & \int^t _0 \bigg(\int_{|x-s_\varepsilon(t)| \leq 3 \varepsilon} |B_\delta|^{\frac{9}{2}} \mathrm{d}x \bigg)^{\frac{2}{3}} \mathrm{d}t\\
\leq & \int^t _0 \bigg(\int_{|x-s_\varepsilon(t)| \leq 3 \varepsilon +\delta} |B|^{\frac{9}{2}} \mathrm{d}x \bigg)^{\frac{2}{3}} \mathrm{d}t \xrightarrow{\delta, \varepsilon \to 0} 0.
\end{split}
\end{equation}

Using $r_\delta,$ we split $V$ as follows.
\begin{equation}\notag
\begin{split}
V= & d_i\iint \big(\nabla \cdot (B \otimes B)_\delta\big) \cdot \big(\nabla \times (A_\delta \phi \chi_\varepsilon)\big)\\
= &d_i\iint r_\delta (B, B) \cdot \nabla \big(\nabla \times (A_\delta \phi \chi_\varepsilon)\big)\\
&+d_i \iint \big((B-B_\delta)\otimes(B-B_\delta)\big) \cdot \nabla \big(\nabla \times (A_\delta \phi \chi_\varepsilon)\big)\\
&-d_i\iint \big(\nabla \cdot (B_\delta\otimes B_\delta)\big) \cdot \big(\nabla \times (A_\delta \phi \chi_\varepsilon)\big)\\
=& V_1 +V_2 -V_3.
\end{split}
\end{equation}

We can prove that $V_1$ and $V_2$ vanish via the same arguments for $IV_1$ and $IV_2.$  
For $V_3,$ we have
\begin{equation}\notag
\begin{split}
V_3 = & d_i \iint \big(\nabla \cdot (B_\delta \otimes B_\delta)\big) \cdot (\nabla \phi \chi_\varepsilon \times A_\delta) \\
&+ d_i \iint \big(\nabla \cdot (B_\delta \otimes B_\delta)\big) \cdot (\phi \nabla \chi_\varepsilon \times A_\delta) \\
=: & V_{31} +V_{32}.
\end{split}
\end{equation}

By standard convergence theorems, as $\delta, \varepsilon \to 0$, $V_{31}$ naturally converges to 
$$\iint \big(\nabla \cdot (B \otimes B) \big) \cdot (\nabla \phi \times A ) = \iint \big( (\nabla \times B) \times B \big) \cdot (\nabla \phi \times A ).$$

As for $V_{32},$ by H\"older's inequality, estimate \eqref{ineqd} and the conditions $A \in L^{\infty} L^6,$ \eqref{acond} and \eqref{bcond}, we have, as $\delta, \varepsilon \to 0,$ that
\begin{equation}\notag
\begin{split}
|V_{32}|\leq & \|A_{\delta}\|_6 \|\nabla B_\delta \phi\chi_\varepsilon\|_{L^{\frac{3}{2}} L^{\frac{18}{5}}} \Bigg(\int^t_0 \bigg(\int_{|x - s_\varepsilon(t)| \leq 3 \varepsilon} |B_\delta|^{\frac{9}{2}} \mathrm{d}x\bigg)^{\frac{2}{3}} \mathrm{d}t \Bigg)^{\frac{1}{3}} \\
\leq & \|A_{\delta}\|_6 \|\nabla B_\delta \phi \chi_\varepsilon \|_{L^{\frac{3}{2}} L^{\frac{18}{5}}} \Bigg(\int^t_0 \bigg(\int_{|x - s_\varepsilon(t)| \leq 3 \varepsilon +\delta } |B|^{\frac{9}{2}} \mathrm{d}x\bigg)^{\frac{2}{3}} \mathrm{d}t \Bigg)^{\frac{1}{3}} \to 0.
\end{split}
\end{equation}

We omit details of the estimates for $VI,$ which are similar to those for $IV,$ while pointing out that $IV_{31}$ is cancelled by its parallel in $VI.$

Integration by parts leads to
\begin{equation}\notag
\begin{split}
VII = & 2\mu \int^t_0 \int_{\R^3} (\nabla A_\delta : \nabla B_\delta) (\phi \chi_\varepsilon)-\mu \int^t_0 \int_{\R^3} (A_\delta \cdot B_\delta) \Delta \phi \chi_\varepsilon\\
& -\mu \int^t_0 \int_{\R^3} (A_\delta \cdot B_\delta) \nabla \phi \nabla \chi_\varepsilon -\mu \int^t_0 \int_{\R^3} (A_\delta \cdot B_\delta) \phi  \Delta \chi_\varepsilon \\
=:& VII_1 + VII_2 + VII_3 + VII_4.
\end{split}
\end{equation}
It's easy to see that as $\delta, \epsilon \to 0,$
$$VII_1 \to 2\mu \iint (\nabla A : \nabla B) \phi \text{ and } VII_2 \to -\mu \iint A\cdot B \Delta\phi.$$

On the other hand, using H\"older's inequality, inequality \eqref{ineqd} and condition \eqref{acond}, we obtain
\begin{equation}\notag
\begin{split}
|VII_3| \leq \mu \varepsilon \sqrt{t}\|A_\delta\|_{L^\infty L^6} \|B_\delta\|_{L^2L^6}\xrightarrow{\delta, \varepsilon \to 0} 0.
\end{split}
\end{equation}
Applying H\"older's inequality, we realize that $VII_4$ vanishes as a consequence of condition \eqref{acond}, estimate \eqref{ineqd} and the fact that $\phi \in \mathscr{D}([0,T]\times \R^3).$
\begin{equation}\notag
\begin{split}
|VII_4| \leq & \mu \|\phi\|_{L^6 L^9} \|A_\delta\|_{L^2 L^\infty}\bigg(\int^t_0 \bigg(\int_{|x-s_\varepsilon(t)| \leq 3 \varepsilon} |B_\delta|^{\frac{9}{2}} \mathrm{d}x \bigg)^{\frac{2}{3}}\mathrm{d}t\bigg)^\frac{1}{3}\\
\leq & \mu \|\phi\|_{L^6 L^9} \|A_\delta\|_{L^2 L^\infty}\bigg(\int^t_0 \bigg(\int_{|x-s_\varepsilon(t)| \leq 3 \varepsilon+\delta} |B|^{\frac{9}{2}} \mathrm{d}x \bigg)^{\frac{2}{3}}\mathrm{d}t\bigg)^\frac{1}{3} \xrightarrow{\delta, \varepsilon \to 0} 0.
\end{split}
\end{equation}

Combining all the estimates above, we recover the generalized helicity equality.
\cbdu

\begin{figure*}[hb]
\begin{tikzpicture}
\draw [<->] (0,5) node[left]{$ \frac{1}{q}$}  --(0,0)  -- (5,0) node[below]{$ \frac{1}{p} $};
\draw  (0.5, 0.5) -- node[below]{$ II $} (0.5, 0.5); 
\draw  (3, 3) -- node[below]{$ I $} (3, 3);
\draw  (4, 0) -- node[below]{$ P_2 $} (4, -0.1); 
\draw  (2, 0) -- node[below]{$ P_4 $} (2, -0.1); 
\draw  (-0.1, 2) -- node[left]{$ P_1 $} (0, 2);
\draw  (-0.1, 1) -- node[left]{$ P_3 $} (0, 1); 
\draw  (0,2)  -- (1.33,1.33) ; \draw   (1.33,1.33)  -- (4,0) ;   
\draw [dashed] (0,1)  -- (2,0) ; 
\draw [->] (4.5,1.5)node[right]{$\frac{2}{q}+\frac{3}{p} = 2$} -- (2.8,0.7)  ;
\fill[pattern=north east lines,opacity=0.5] (0,1)-- (0,2) --(4,0)-- (2,0) ;
\draw[->]    (1,-0.5)node[left]{ $\frac{2}{q}+ \frac{3}{p}=1$ } --(1.5,0.3) ;
\end{tikzpicture}
\caption{Conditions on parameters $p$ and $q$ indicating rigidity.}
\label{fig-para}
\end{figure*}
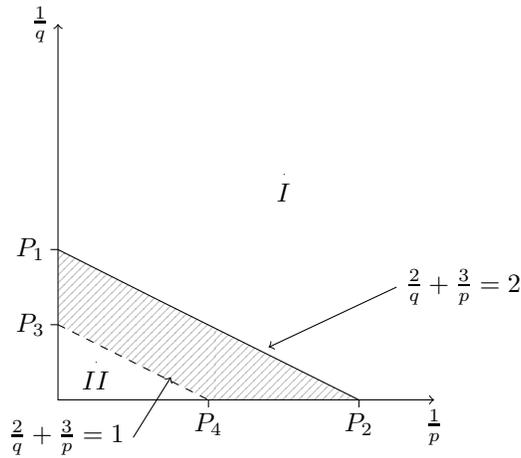

\bigskip

\section{Discussions}
In this section, we further analyze the main results. We can visualize the condition in Theorem \ref{thm1} that guarantees uniqueness in Figure \ref{fig-para}.  Theorem \ref{thm1} says that if a Leray-Hopf weak solution $B(t)$ satisfies  $\nabla\times B\in L^q(B^{-1+\frac2q+\frac3p}_{p,\infty})$ for $(\frac1p, \frac1q)$ in the shaded region of Figure \ref{fig-para}, it is unique in this class. In this figure, $P_1$, $P_2$, $P_3$, and $P_4$ correspond to the spaces $L^{3/2}(B^1_{\infty,\infty})$, $L^\infty(B^1_{1,\infty})$, $L^3(B^0_{\infty,\infty})$, and $L^\infty(B^0_{2,\infty})$ respectively. Indicated by Remark 1.4, if $\nabla B$ belongs to the Lebesgue space counterpart $L^q(L^p)$ with $(\frac1p, \frac1q)$ in region II, then $B$ is regular and hence unique in this class. While in region I, the uniqueness or lack of uniqueness remains open.

Implied by Theorem \ref{consh}, if a finite energy solution $B(t)$ satisfies $\nabla B\in L^3(B^{-2/3}_{3,\infty})$, then it conserves the magnetic helicity. The corresponding point of exponent parameter falls in the region I. Thus, one can see that the conservation of magnetic helicity represents weaker rigidity than that of uniqueness, and uniqueness represents a weaker rigidity than that of regularity.


\section{Appendix}

\subsection{Littlewood-Paley theory}

Here we give a concise review of Littlewood-Paley theory. For a complete description of the theory and its applications, readers are referred to the books \cite{BCD} and \cite{G}.

We construct a family of smooth functions $\{\varphi_q \}_{q=-1}^\infty$ with annular support that forms a dyadic partition of unity in the frequency space, defined as 
\begin{equation}\notag
\varphi_q(\xi)=
\begin{cases}
\varphi(\lambda_q^{-1}\xi)  \ \ \ \mbox { for } q\geq 0,\\
\chi(\xi) \ \ \ \mbox { for } q=-1,
\end{cases}
\end{equation}
where  $\lambda_q = 2^q,$ $\varphi(\xi)=\chi(\xi/2)-\chi(\xi)$ and $\chi\in C_0^\infty(\R^n)$ is a nonnegative radial function chosen in a way such that 
\begin{equation}\notag
\chi(\xi)=
\begin{cases}
1, \ \ \mbox { for } |\xi|\leq\frac{3}{4}\\
0, \ \ \mbox { for } |\xi|\geq 1.
\end{cases}
\end{equation}

Introducing the functions $\tilde h:=\mathcal{F}^{-1}(\chi)$ and $h:= \mathcal{F}^{-1}(\varphi),$ we define the inhomogeneous Littlewood-Paley projections for $u \in \mathscr{S}^{'}(\R^n)$ as
\begin{equation}\notag
u_q:=\Delta_qu= \mathcal{F}^{-1} (\varphi_q(\xi)\hat u(\xi) ) =  \begin{cases}
\displaystyle \lambda_q^n \int_{\R^3} h(\lambda_q y) u(x-y) \mathrm{d}y, \  q \geq 0, \\
\displaystyle\int_{\R^3} \tilde h(y) u(x-y) \mathrm{d}y, \ q =-1.
\end{cases}
\end{equation}

Formally, the identity
$$
u=\sum_{q=-1}^\infty u_q
$$
holds at least in the sense of distributions. To simplify the notation, we denote
$$u_{\leq Q}=\sum_{q=-1}^Qu_q.$$

We recall Bernstein's inequality, whose proof can be found in \cite{BCD}.
\begin{Lemma}\label{brnst}
Let $n$ be the space dimension and $1\leq s\leq r\leq \infty$. Then for all tempered distributions $u$, 
$$
\|u_q\|_{r}\lesssim \lambda_q^{n(\frac{1}{s}-\frac{1}{r})}\|u_q\|_{s}.
$$
\end{Lemma}

\subsection{Proof of Theorem \ref{consh}}

We give a proof of the positive side of the analogue of Onsager's conjecture for the non-resistive electron-MHD system, written as follows.
\begin{equation}\label{AEMHD}
\begin{cases}
A_t= d_i (\nabla \times B)\times B,\\
\nabla \times A = B, \ \nabla \cdot B =0, \ t \in \R^+, x \in \R^3 (\text{ or } \mathbb T^3).
\end{cases}
\end{equation}
Our proof follows that in \cite{CCFS}, where the positive side of Onsager's conjecture was confirmed by the result that any weak solution $u \in L^3(0,T; B^{1/3}_{3, c(\mathbb{N})})$ to 3D Euler's equations conserves energy.

Clearly, for regular solutions to \eqref{AEMHD}, energy and magnetic helicity are conserved.
\begin{equation}\notag
\frac{\mathrm{d}}{\mathrm{d}t}\mathcal{H}(t)=\frac{\mathrm{d}}{\mathrm{d}t}\int_{\R^3} \big(A \cdot B\big)(x,t) \mathrm{d}x =0; \ \ 
\frac{\mathrm{d}}{\mathrm{d}t}\mathcal{E}(t)=\frac{1}{2}\frac{\mathrm{d}}{\mathrm{d}t}\int_{\R^3} |B(x,t)|^2 \mathrm{d}x =0.
\end{equation}

On the other hand, to our knowledge, the existence of weak solutions to system \eqref{AEMHD} remains an open question at this time. We say that $(A, B)$ is a weak solution to system \eqref{AEMHD}, if $(A, B)$ is a pair of divergence-free vector fields satisfying the equations in the sense of distributions. 

To this end, we shall show that the total helicity flux of any divergence-free vector field $B \in B^{1/3}_{3, c(\mathbb{N})}$ vanishes, which in turn implies conservation of helicity. 

To start, we define the truncated helicity flux as
\begin{equation}
\mathcal{H}_Q =2\int_{\R^3} ((\nabla \times B) \times B)_{\leq Q}\cdot   B_{\leq Q} \ \mathrm{d}x
\end{equation}
and the truncated energy flux as 
\begin{equation}
\Pi_Q =\int_{\R^3} ((\nabla \times B) \times B)_{\leq Q}\cdot (\nabla \times  B_{\leq Q})\ \mathrm{d}x.
\end{equation}
We note that  
$$
\mathcal{H}_Q(t)=\frac{\mathrm{d}}{\mathrm{d}t}\int_{\R^3} A_{\leq Q} \cdot B_{\leq Q} \ \mathrm{d}x \text{ and } \Pi_Q(t)=\frac{1}{2}\frac{\mathrm{d}}{\mathrm{d}t}\|B_{\leq Q}(t)\|_2^2,
$$
provided that $B$ is a weak solution to \eqref{AEMHD}.

We introduce the localization kernels 
\begin{equation}
\mathcal{K}(q)=\begin{cases}
\lambda_q^{2/3}, \ q \leq 0,\\
\lambda_q^{-4/3}, \ q > 0;
\end{cases}
\end{equation}
and
\begin{equation}
\kappa(q)=\begin{cases}
\lambda_q^{4/3}, \ q \leq 0,\\
\lambda_q^{-2/3}, \ q > 0.
\end{cases}
\end{equation}

For $B \in \mathcal{S}^{'},$ we define
\begin{equation}
b_q:= \lambda_q^{1/3}\|B_q\|_3 \text{ and } \beta_q:= \lambda_q^{2/3}\|B_q\|_3,
\end{equation}
and denote the sequences $\{b_q^2\}_{q = -1}^\infty$ and $\{\beta_q^2\}_{q=-1}^\infty$ by $b^2$ and $\beta^2,$ respectively.

By vector identities and integration by parts, we have
\begin{equation}\notag
\begin{split}
\mathcal{H}_Q = & 2\int_{\R^3} ((\nabla \times B) \times B)_{\leq Q}\cdot   B_{\leq Q} \ \mathrm{d}x\\
= & 2\int_{\R^3} \Big(\nabla \cdot (B \otimes B) - \frac{1}{2}\nabla|B|^2\Big)_{\leq Q}\cdot   B_{\leq Q} \ \mathrm{d}x\\
= & -2\int_{\R^3}  (B \otimes B)_{\leq Q} :  \nabla B_{\leq Q} \ \mathrm{d}x +\int_{\R^3}|B|^2_{\leq Q}\nabla \cdot B_{\leq Q} \ \mathrm{d}x\\
=: &\ 2\mathcal{H}_Q^1 + \mathcal{H}_Q^2.
\end{split}
\end{equation}

We shall estimate $\mathcal{H}_Q^1$ only, as $\mathcal{H}_Q^2$ can be estimated in a similar way. 

Introducing a bilinear form
$$r_Q (B, B): = \int_{\R^3} \tilde h(y) (B(x-y)-B(x))\otimes (B(x-y)-B(x)) \mathrm{d}y,$$
we can split $(u \otimes u)_{\leq Q}$ into three parts.
\begin{equation}\notag
\begin{split}
(u \otimes u)_{\leq Q} = & r_Q (B, B) - (B - B_{\leq Q}) \otimes (B - B_{\leq Q}) + B_{\leq Q}\otimes B_{\leq Q}.
\end{split}
\end{equation}

Integration by parts yields
\begin{equation}\notag
\begin{split}
\mathcal{H}_Q^1 = & \int_{\R^3} r_Q (B, B) \cdot \nabla B_{\leq Q} \mathrm{d}x\\
& - \int_{\R^3} (B - B_{\leq Q}) \otimes (B - B_{\leq Q})\cdot \nabla B_{\leq Q} \mathrm{d}x \\
= & \mathcal{H}_Q^{11} + \mathcal{H}_Q^{12}.
\end{split}
\end{equation}

By H\"older's inequality, we have
\begin{equation}\notag
\mathcal{H}_Q^{11} \leq  \|r_Q (B, B)\|_{\frac{3}{2}}\|\nabla B_{\leq Q}\|_3,
\end{equation}
and
\begin{equation}\notag
\begin{split}
\|r_Q (B, B)\|_{\frac{3}{2}} \leq & \int_{\R^3}\big|\tilde h_Q(y)\big|\|B(\cdot-y)-B(\cdot)\|_3^2 \mathrm{d}y.
\end{split}
\end{equation}

Separating the lower and higher frequencies, we obtain
\begin{equation}\notag
\begin{split}
\|B(\cdot-y)-B(\cdot)\|_3^2 \leq & \Big(\sum_{q \leq Q}|y|^2 \lambda_q^2\|B_q\|_3^2  +\sum_{q >Q}\|B_q\|_3^2\Big) \\
\leq & \lambda_Q^{\frac{4}{3}}|y|^2 \sum_{q \leq Q}\lambda_{Q-q}^{-\frac{4}{3}}b_q^2 + \lambda^{-\frac{2}{3}}_Q \sum_{q >Q} \lambda_{Q-q}^{\frac{2}{3}}b^2_q\\
\leq & \Big(\lambda_Q^{\frac{4}{3}}|y|^2 + \lambda^{-\frac{2}{3}}_Q\Big)(K \ast b^2)(Q).
\end{split}
\end{equation}

It follows that
\begin{equation}\notag
\begin{split}
|\mathcal{H}_Q^{11}| \leq & (K \ast b^2)(Q) \Big( \int_{\R^3}\big|\tilde h_Q(y)\big|\lambda_Q^{\frac{4}{3}}|y|^2 \mathrm{d}y + \lambda^{-\frac{2}{3}}_Q \Big) \|\nabla B_{\leq Q}\|_3 \\
\leq  & (K \ast b^2)(Q) \Big( \int_{\R^3} \big|\tilde h_Q(y)\big|\lambda_Q^{\frac{4}{3}}|y|^2 \mathrm{d}y + \lambda^{-\frac{2}{3}}_Q \Big) \Big( \sum_{q \leq Q} \lambda_q^2 \|B_q\|_3^2 \Big)^{\frac{1}{2}} \\
\leq  & (K \ast b^2)(Q) \Big( \int_{\R^3} \big|\tilde h_Q(y)\big|\lambda_Q^{\frac{4}{3}}|y|^2 \mathrm{d}y + \lambda^{-\frac{2}{3}}_Q \Big) \Big( \sum_{q \leq Q} \lambda_q^\frac{4}{3} b_q^2 \Big)^{\frac{1}{2}} \\
\leq  & (K \ast b^2)(Q) \lambda^{-\frac{2}{3}}_Q \Big( \sum_{q \leq Q}\lambda_q^\frac{4}{3} b_q^2 \Big)^{\frac{1}{2}} \\
\leq  & (K \ast b^2)^\frac{3}{2}(Q).
\end{split}
\end{equation}
As $B \in L^3(0,T; B^{2/3}_{3, c(\mathbb{N})}),$ it is clear that $$\lim_{Q \to \infty}|\mathcal{H}_Q^{11}|=0. $$

Analogously, we have 
\begin{equation}\notag
\begin{split}
|\mathcal{H}_Q^{12}| \leq & \|B-B_{\leq Q}\|_3^2 \|\nabla B_{\leq Q}\|_3\\
\leq & \Big( \sum_{q > Q}\|B_q\|_3^2 \Big)\Big(\sum_{q \leq Q}\lambda^2_q \|B_q\|_3^2 \Big)^\frac{1}{2}\\
\leq  & (K \ast b^2)^\frac{3}{2}(Q),
\end{split}
\end{equation}
indicating that $|\mathcal{H}_Q^{12}|$ vanishes.

\begin{Remark} We can prove the following theorem regarding energy conservation for weak solutions to system \eqref{AEMHD} via the same approach as above.
\begin{Theorem}\label{thm6}
Let $B \in L^3(0,T; B^{2/3}_{3, c(\mathbb{N})})\cap C_w(0,T; L^2)$ be a weak solution to \eqref{AEMHD}, then $B$ conserves the magnetic energy $\mathcal{E}.$
\end{Theorem}
\end{Remark}

\end{document}